\documentclass[11pt,twoside,reqno,centertags,draft]{amsart}
\usepackage{amsfonts}
\usepackage{color,enumitem,graphicx}
\usepackage[colorlinks=true,urlcolor=blue,
citecolor=red,linkcolor=blue,linktocpage,pdfpagelabels,
bookmarksnumbered,bookmarksopen]{hyperref}

  \usepackage{amsmath,amsthm,amsfonts,amssymb}

\begin{document}
\title{ Nonlinear Choquard equations involving\\
\smallskip
 nonlocal operators}
\date{}
\maketitle

\vspace{ -1\baselineskip}

{\small
\begin{center}

\medskip

  {\sc  Wanwan Wang}
  \medskip

Department of Mathematics, Jiangxi Normal University,\\
 Nanchang, Jiangxi 330022, PR China
\end{center}
}

\renewcommand{\thefootnote}{}

\footnote{E-mail address: wwwang2014@yeah.net (W. Wang).}

\begin{quote}
{\bf Abstract.} In this paper, we study  nonlinear Choquard equations
\begin{equation}\label{eq 1a1-}
(-\Delta+id)^{\frac{1}{2}}u=(I_\alpha*{|u|^p})|u|^{p-2}u\ \ {\rm in} \ \ \mathbb{R}^N, \ \ \ u\in H^{\frac{1}{2}}(\mathbb{R}^N),
\end{equation}
where $(-\Delta+id)^\frac{1}{2}$ is a nonlocal operator, $p>0$, $N\geq2$ and $I_\alpha$ is the Riesz potential with  order $\alpha\in(0,N)$.
We show that there is a ground state solution to problem (\ref{eq 1a1-}) if $\frac{N+\alpha}{N}<p<\frac{N+\alpha}{N-1}$ and
 no solution to problem (\ref{eq 1a1-}) if $0<p\leq\frac{N+\alpha}{N+1}$ or $p\geq\frac{N+\alpha}{N-1}$. Furthermore, the existence of infinity many solutions to
 problem (\ref{eq 1a1-}) is discussed when $p$ satisfies that $\frac{N+\alpha}{N}<p<\frac{N+\alpha}{N-1}$.

\end{quote}

\medskip

\newcommand{\N}{\mathbb{N}}
\newcommand{\R}{\mathbb{R}}
\newcommand{\Z}{\mathbb{Z}}

\newcommand{\cA}{{\mathcal A}}
\newcommand{\cB}{{\mathcal B}}
\newcommand{\cC}{{\mathcal C}}
\newcommand{\cD}{{\mathcal D}}
\newcommand{\cE}{{\mathcal E}}
\newcommand{\cF}{{\mathcal F}}
\newcommand{\cG}{{\mathcal G}}
\newcommand{\cH}{{\mathcal H}}
\newcommand{\cI}{{\mathcal I}}
\newcommand{\cJ}{{\mathcal J}}
\newcommand{\cK}{{\mathcal K}}
\newcommand{\cL}{{\mathcal L}}
\newcommand{\cM}{{\mathcal M}}
\newcommand{\cN}{{\mathcal N}}
\newcommand{\cO}{{\mathcal O}}
\newcommand{\cP}{{\mathcal P}}
\newcommand{\cQ}{{\mathcal Q}}
\newcommand{\cR}{{\mathcal R}}
\newcommand{\cS}{{\mathcal S}}
\newcommand{\cT}{{\mathcal T}}
\newcommand{\cU}{{\mathcal U}}
\newcommand{\cV}{{\mathcal V}}
\newcommand{\cW}{{\mathcal W}}
\newcommand{\cX}{{\mathcal X}}
\newcommand{\cY}{{\mathcal Y}}
\newcommand{\cZ}{{\mathcal Z}}

\newcommand{\abs}[1]{\lvert#1\rvert}
\newcommand{\xabs}[1]{\left\lvert#1\right\rvert}
\newcommand{\norm}[1]{\lVert#1\rVert}

\newcommand{\loc}{\mathrm{loc}}
\newcommand{\p}{\partial}
\newcommand{\h}{\hskip 5mm}
\newcommand{\ti}{\widetilde}
\newcommand{\D}{\Delta}
\newcommand{\e}{\epsilon}
\newcommand{\bs}{\backslash}
\newcommand{\ep}{\emptyset}
\newcommand{\su}{\subset}
\newcommand{\ds}{\displaystyle}
\newcommand{\ld}{\lambda}
\newcommand{\vp}{\varphi}
\newcommand{\wpp}{W_0^{1,\ p}(\Omega)}
\newcommand{\ino}{\int_\Omega}
\newcommand{\bo}{\overline{\Omega}}
\newcommand{\ccc}{\cC_0^1(\bo)}
\newcommand{\iii}{\opint_{D_1}D_i}

\numberwithin{equation}{section}

\vskip 0.2cm \arraycolsep1.5pt
\newtheorem{lemma}{Lemma}[section]
\newtheorem{theorem}{Theorem}[section]
\newtheorem{definition}{Definition}[section]
\newtheorem{proposition}{Proposition}[section]
\newtheorem{remark}{Remark}[section]
\newtheorem{corollary}{Corollary}[section]

\setcounter{equation}{0}
\section{Introduction}

Our purpose of this paper is to consider the
 solutions of nonlinear Choquard equations
\begin{equation}\label{eq 1}
(-\Delta+id)^{\frac{1}{2}}u=(I_\alpha*{|u|^p})|u|^{p-2}u\ \ {\rm in} \ \ \mathbb{R}^N,\ \ \ u\in H^{\frac{1}{2}}(\mathbb{R}^N),
\end{equation}
where  $p>0$, $N\ge2$, $I_\alpha:\mathbb{R}^N \backslash\{0\}\to \R$ is
the Riesz potential with  order $\alpha\in(0,N)$ given  by
$$I_\alpha(x)=\frac{\Gamma(\frac{N-\alpha}{2})}{2^\alpha \pi^\frac{N}{2}\Gamma(\frac{\alpha}{2})|x|^{N-\alpha}},$$
here $\Gamma$ is gamma function, see \cite{RM}.
The nonlocal operator  $(-\Delta+id)^{\frac{1}{2}}$ can be characterized as $\mathcal{F}((-\Delta+id)^{\frac{1}{2}}\phi)(\xi)=(1+|\xi|^2)^{\frac12}\mathcal{F}(\phi)(\xi)$,
here $\mathcal{F}$ is the Fourier transform.
The Hilbert space $H^{\frac{1}{2}}(\mathbb{R}^N)$ is defined as
$$H^{\frac{1}{2}}(\mathbb{R}^N)=\{u\in L^2(\mathbb{R}^N): \ (-\Delta+id)^\frac{1}{4}u\in L^2(\mathbb{R}^N)\}$$
with the norm
$$\|u\|_{H^{\frac{1}{2}}(\mathbb{R}^N)}=\|(-\Delta+id)^\frac{1}{4}u\|_{L^2(\mathbb{R}^N)}.$$

As early as in 1954, in a pioneering  work of Pekar \cite{P} where described the quantum mechanics of a polaron, the nonlinear Choquard equation
$$
\left\{ \arraycolsep=1pt
\begin{array}{lll}
 -\Delta u+u=(I_2\ast|u|^2)u\ \ \ &
{\rm in} \quad \R^3,\\[2mm]
\lim_{|x|\to+\infty}u(x)=0
\end{array}
\right.
$$
is appeared. For general case, Moroz-Van Schaftingen in \cite{VJ} studied the problem
$$
\left\{ \arraycolsep=1pt
\begin{array}{lll}
 -\Delta u+u=(I_\alpha\ast|u|^p)|u|^{p-2}u\ \ \ &
{\rm in} \quad \R^N,\\[2mm]
\lim_{|x|\to+\infty}u(x)=0,
\end{array}
\right.
$$
they obtained the results of existence, qualitative properties and decay asymptotic.
In this paper, we consider the related nonlocal problem (\ref{eq 1}).
To state our results, we first introduce
Hardy-Littlewood-Sobolev inequality which states that if $s\in(1,\frac{N}{\alpha})$,
 then for every $v\in L^s(\mathbb{R}^N)$, $I_\alpha*v\in L^\frac{Ns}{N-\alpha s}(\mathbb{R}^N)$ and
%\begin{eqnarray*}\label{10e}
$$\int_{\mathbb{R}^N}|I_\alpha*v|^{\frac{Ns}{N-\alpha s}}dx\leq C(\int_{\mathbb{R}^N}|v|^sdx)^\frac{N}{N-\alpha s},$$
%\end{eqnarray*}
where $C>0$ depends  on $\alpha$, $N$ and $s$. We know that the fractional Sobolev embedding $H^{\frac{1}{2}}(\mathbb{R}^N)\subset L^t(\mathbb{R}^N)$ for $t\in[2,2^\sharp]$, where $2^\sharp:=\frac{2N}{N-1}$, also note that $H^{\frac{1}{2}}(\mathbb{R}^N)\subset L^\frac{2Np}{N+\alpha}(\mathbb{R}^N)$ if and only if $\frac{N+\alpha}{N}<p<\frac{N+\alpha}{N-1}$.
Now we state our main theorem.

\begin{theorem}\label{th1}
Assume that $N\ge2$ and $\alpha\in(0,N)$.

$(i)$ If $\frac{N+\alpha}{N}<p<\frac{N+\alpha}{N-1}$,
then there exists a $C^2$ positive ground state solution to problem (\ref{eq 1}).

$(ii)$ If $0<p\leq\frac{N+\alpha}{N+1}$ or $p\geq\frac{N+\alpha}{N-1}$, then there is no nontrivial solution to problem (\ref{eq 1}).
\end{theorem}

To prove the existence of solutions in Theorem \ref{th1} when $\frac{N+\alpha}{N}<p<\frac{N+\alpha}{N-1}$, we apply the critical points theory to the associated minimizing problem
\begin{equation}\label{4e}
M_p=\inf\left\{\int_{\mathbb R^N}\left|(-\Delta+id)^\frac{1}{4}u\right|^2dx: \int_{\mathbb{R}^N}(I_\alpha*|u|^p)(x)|u(x)|^p\,dx=1\right\}.
\end{equation}
We note that the minimization of $M_p$ is a  nontrivial solution of
problem (\ref{eq 1}). Here we use the  concentration compactness argument and a nonlocal version of Brezis-Lieb lemma to prove that $M_p$ can be achieved.
Then we establish Pohoz\v{a}ev identity to obtain the nonexistence results in Theorem \ref{th1}.

\begin{theorem}\label{th2}
Let  $N\ge2$,  $\alpha\in(0,N)$ and $\frac{N+\alpha}{N}<p<\frac{N+\alpha}{N-1}$.
Then there exists infinitely many distinct solutions to problem (\ref{eq 1}).
\end{theorem}

\setcounter{equation}{0}
\section{Preliminaries}
In this section, we introduce some lemmas.

\begin{lemma}\label{lemma 2} \cite {W}
Let $\Omega $ be a domain in $\mathbb{R}^N$, $s>1$ and $\{w_m\}_{m\in \mathbb{N}}$ be a bounded sequence in $L^s(\Omega)$.
If $w_m\to w$ almost everywhere on $\Omega$ as $m\to \infty$, then for every $q\in[1,s]$, we have that
$$\lim_{m\to\infty}{\int_\Omega|{|w_m|^q-|w_m-w|^q-|w|^q}|^\frac{s}{q}}\,dx=0.$$
 \end{lemma}

\begin{lemma}\label{lemma 1}
Let  $\alpha \in(0,N)$, $\frac{N+\alpha}{N}<p<\frac{N+\alpha}{N-1}$ and $\{w_m\}_{m\in \mathbb{N}}$ be a bounded sequence in
$L^{\frac{2Np}{N+\alpha}}(\mathbb{R}^N)$. Assume that

\quad (i) ${w_m}$ weakly converges to ${w}$ in $L^\frac{2Np}{N+\alpha}(\mathbb{R}^N)$;

\quad (ii)  ${w_m}\to w$ almost everywhere on  $\mathbb{R}^N$.\\
Then
\begin{eqnarray*}
&&\lim_{m\to \infty}\left[\int_{\mathbb{R}^N}(I_\alpha*|w_m|^p)(x)|w_m(x)|^p\,dx-\int _{\mathbb{R}^N}(I_\alpha*|w_m-w|^p)(x)|(w_m-w)(x)|^p\,dx\right]
\\&&=\int_{\mathbb{R}^N}(I_\alpha*|w|^p)(x)|w(x)|^p\,dx.
\end{eqnarray*}
 \end{lemma}
\noindent{\bf Proof.}
We observe that
 \begin{eqnarray*}
 && \int_{\mathbb{R}^N}(I_\alpha*|w_m|^p)(x)|w_m(x)|^p\,dx-\int _{\mathbb{R}^N}(I_\alpha*|w_m-w|^p)(x)|(w_m-w)(x)|^p\,dx
\\&=&\int_{\mathbb{R}^N}(I_\alpha*(|w_m|^p-|w_m-w|^p))(x)(|w_m(x)|^p-|(w_m-w)(x)|^p)\,dx
\\&&+2\int_{\mathbb{R}^N}(I_\alpha*(|w_m|^p-|w_m-w|^p))(x)|(w_m-w)(x)|^p\,dx.
\end{eqnarray*}
By the H\"{o}lder inequality, we have that
\begin{eqnarray*}
 && \int_{\mathbb{R}^N}(I_\alpha*(|w_m|^p-|w_m-w|^p))(x)|(w_m-w)(x)|^p\,dx
\\&=&\int_{\mathbb{R}^N}(I_\alpha*(|w_m|^p-|w_m-w|^p-|w|^p))(x)|(w_m-w)(x)|^p\,dx
\\&&+\int_{\mathbb{R}^N}(I_\alpha*|w|^p)(x)|(w_m-w)(x)|^p\,dx
\\&\leq &\left({\int_{\mathbb{R}^N}|\left(I_\alpha*(|w_m|^p-|w_m-w|^p-|w|^p)\right)|^\frac{2N}{N-\alpha}}(x)\,dx\right)^\frac{N-\alpha}{2N}
\\&&\left(\int_{\mathbb{R}^N}(|(w_m-w)(x)|^p)^\frac{2N}{N+\alpha}\,dx\right)^\frac{N+\alpha}{2N}+\int_{\mathbb{R}^N}(I_\alpha*|w|^p)(x)|(w_m-w)(x)|^p\,dx.
\end{eqnarray*}
Using Lemma \ref{lemma 2} with $q=p$ and $s=\frac{2Np}{N+\alpha}$, we know that $|w_m|^p-|w_m-w|^p\to |w|^p$,
strongly in $L^\frac{2N}{N+\alpha}(\mathbb{R}^N)$ as $m\to\infty$. By the Hardy-Littlewood-Sobolev inequality, this implies that
$I_\alpha*(|w_m|^p-|w_m-w|^p)\to I_\alpha*|w|^p$ in $L^\frac{2N}{N-\alpha}(\mathbb{R}^N)$ as $m\to\infty$. Since $|w_m-w|^p\rightharpoonup 0$
in $L^\frac{2N}{N+\alpha}(\mathbb{R}^N)$ as $m\to\infty$, then $\int_{\mathbb{R}^N}(I_\alpha*(|w_m|^p-|(w_m-w)|^p))(x)|(w_m-w)(x)|^p\,dx\to0$
as $m\to\infty$. This ends the proof. \hfill$\Box$

\setcounter{equation}{0}
\section{Ground state solution}

In this section, we study the existence of ground state solutions of problem (\ref{eq 1}). To this end,
let us consider the following minimizing problem
\begin{equation}\label{2c}
M_p=\inf\left\{\int_{\mathbb R^N}\left|(-\Delta+id)^\frac{1}{4}u\right|^2dx: \int_{\mathbb{R}^N}(I_\alpha*|u|^p)(x)|u(x)|^p\,dx=1\right\}.
\end{equation}
\begin{proposition}\label{proposition 3.1}
The minimizing problem $M_p$ is achieved by a function $v\in {H^{\frac{1}{2}}(\mathbb{R}^N)}$, which is a solution of problem (\ref{eq 1})
up to a translation.
\end{proposition}
To prove this result, we first introduce some lemma as follows.
\begin{lemma}\label{lemma 2.1}
Let $2\leq q<2^\sharp=\frac{2N}{N-1}$ and $r>0$. Suppose that $\{v_m\}_{m\in\mathbb{N}}$ is a bounded sequence in $H^{\frac{1}{2}}(\mathbb{R}^N)$
and
$$\sup_{z\in{\mathbb{R}^N}}\int _{B_r(z)}|v_m(x)|^q\,dx\to 0$$
as $m\to\infty$.
Then
 $$\int_{\mathbb{R}^N}(I_\alpha*|v_m|^p)(x)|v_m(x)|^p\,dx\to0$$
as $m\to\infty$ for $\frac{(N+\alpha)q}{2N}<p<\frac{N+\alpha}{N-1}$.
\end{lemma}

\noindent {\bf Proof.}
Let $l=t=\frac{2N}{N+\alpha}$, by Hardy-Littlewood-Sobolev inequality, we have that
\begin{eqnarray*}
 && \int_{\mathbb{R}^N}(I_\alpha*|v_m|^p)(x)|v_m(x)|^p\,dx\ =\int_{\mathbb{R}^N}\int_{\mathbb{R}^N}\frac{|v_m(x)|^p|v_m(z)|^p}{|x-z|^{N-\alpha}}\,dxdz
\\& \leq& c_1\|{|v_m|^p}\|_{L^l(\mathbb{R}^N)}\|{|v_m|^p}\|_{L^t(\mathbb{R}^N)} \  =c_1\left(\int_{\mathbb{R}^N}|v_m(x)|^\frac{2Np}{N+\alpha}\,dx\right)^\frac{N+\alpha}{N}
\end{eqnarray*}
and
\begin{eqnarray*}
\left(\int_{B_r(z)}|v_m(x)|^\frac{2Np}{N+\alpha}\,dx\right)^\frac{N+\alpha}{2Np}
&\leq & \left[\left(\int_{B_r(z)}|v_m(x)|^q\,dx\right)^\frac{1}{q}\right]^{1-\lambda}\left[\left(\int_{B_r(z)}|v_m(x)|^{2^\sharp}\,dx\right)^\frac{1}{2^\sharp}\right]^\lambda
\\&\leq & c_2{\|v_m\|}^{1-\lambda}_{L^q(B_r(z))}\left(\int_{\mathbb{R}^N}|(-\Delta+id)^\frac{1}{4}v_m|^2\,dx\right)^\frac{\lambda}{2},
\end{eqnarray*}
where $c_1,c_2>0$ and
$$\lambda=\frac{\frac{2Np}{N+\alpha}-q}{\frac{2N}{N-1}-q}\cdot\frac{\frac{2N}{N-1}}{\frac{2Np}{N+\alpha}}.$$
 Choosing some suitable $q$ and $p$ such that $\lambda=\frac{2}{\frac{2Np}{N+\alpha}}$,
we obtain that
\begin{eqnarray*}
\int_{B_r(z)}|v_m(x)|^\frac{2Np}{N+\alpha}\,dx \leq  c_2^\frac{2Np}{N+\alpha}{\|v_m\|}^{(1-\lambda)\frac{2Np}{N+\alpha}}_{L^q(B_r(z))}\left(\int_{\mathbb{R}^N}|(-\Delta+id)^\frac{1}{4}v_m|^2\,dx\right).
\end{eqnarray*}
Now, covering $\mathbb{R}^N$ by balls of radius $r$, in such a way that each point of $\mathbb{R}^N$ is contained in at most $N+1$ balls, we have that
\begin{eqnarray*}
&&\int_{B_r(z)}|v_m(x)|^\frac{2Np}{N+\alpha}\,dx
\\&\leq & (N+1)c_2^\frac{2Np}{N+\alpha}\sup_{y\in{\mathbb{R}^N}}[\int_{B_r(z)}{|v_m(x)|}^q\,dx]^\frac{(1-\lambda)\frac{2Np}{N+\alpha}}{q}\left(\int_{\mathbb{R}^N}
|(-\Delta+id)^\frac{1}{4}v_m|^2\,dx\right).
\end{eqnarray*}
Then
\begin{eqnarray*}
&&\int_{\mathbb{R}^N}(I_\alpha*|v_m|^p(x))|v_m(x)|^p\,dx
\\&\leq & c_3 \left((\sup_{z\in{\mathbb{R}^N}}\int_{B_r(z)}|v_m(x)|^q\,dx)^{\frac{1}{q}(\frac{2Np}{N+\alpha}-2)}
(\int_{\mathbb{R}^N}|(-\Delta+id)^\frac{1}{4}v_m|^2\,dx)\right)^\frac{N+\alpha}{N}.
\end{eqnarray*}
The proof is complete.\hfill$\Box$

We now  prove proposition \ref{proposition 3.1}.

\noindent{\bf Proof of Proposition \ref{proposition 3.1}}
\ Let $\{v_m\}_{m\in\mathbb{N}}\subset H^{\frac{1}{2}}(\mathbb{R}^N)$ be a minimizing sequence of $M_p$, which satisfies that
$$\int_{\mathbb{R}^N}(|\nabla v_m(x)|^2+|v_m(x)|^2)\,dx \to M_p$$
and
$$\int_{\mathbb{R}^N}(I_\alpha*|v_m|^p)(x)|v_m(x)|^p\,dx=1.$$
By Lemma \ref{lemma 2.1}, there exists $\delta>0$ such that
$$\delta=\lim_{m\to\infty}\sup_{\tilde{x}\in\mathbb{R}^N}\int_{B_1(\tilde{x})}|v_m(x)|^\frac{2Np}{N+\alpha}\,dx>0.$$
So we may find $x_m\in\mathbb{R}^N$ such that
$$\int_{B_1(x_m)}|v_m(x)|^\frac{2Np}{N+\alpha}\,dx>\frac{\delta}{2}.$$
Let $w_m(x)=v_m(x-x_m)$, we have that
$$\int_{\mathbb{R}^N}(I_\alpha*|w_m|^p)(x)|w_m(x)|^p\,dx=1,\ \ \ \int_{\mathbb R^N}|(-\Delta+id)^\frac{1}{4}w_m|^2\,dx \to M_p,$$
which yields, up to a subsequence, that  $w_m\rightharpoonup w$ in $H^{\frac12}(\R^N)$ and $w_m \to w$ almost everywhere on $\mathbb{R}^N$.
Then
$$\int_{B_1(0)}|w_m(x)|^\frac{2Np}{N+\alpha}\,dx>\frac{\delta}{2}.$$
Since $H_{loc}^\frac{1}{2}(\mathbb{R}^N)\subset L_{loc}^\frac{2Np}{N+\alpha}(\mathbb{R}^N)$ is compact, so we have that$\int_{B_1(0)}|w(x)|^\frac{2Np}{N+\alpha}\,dx\geq\frac{\delta}{2}>0$,we can claim that $w\neq 0$ almost everywhere on $\mathbb{R}^N$.
Then $\int_{\mathbb{R}^N}(I_\alpha*|w|^p)(x)|w(x)|^p\,dx\neq 0$ almost everywhere on $\mathbb{R}^N$.

%We next claim that $\int_{\mathbb{R}^N}(I_\alpha*|w|^p)(x)|w(x)|^p\,dx\neq 0$. In fact, we observe that
%\begin{eqnarray*}
%&&\int_{B_1(x_m)}(I_\alpha*|w_m|^p)(x)|w_m(x)|^p\,dx
%=\int_{B_1(x_m)}(I_\alpha*|v_m|^p)(x-x_m)|v_m(x-x_m)|^p\,dx
%\\&=&\int_{B_1(0)}(I_\alpha*|v_m|^p)(x)|v_m(x)|^p\,dx
%=\int_{B_1(0)}\int_{\mathbb{R}^N}\frac{|v_m(x)|^p}{|x-z|^{N-\alpha}}|v_m(x)|^p\,dzdx
%\\&=&\int_{\mathbb{R}^N}\frac{|v_m(x)|^p}{|x-z|^{N-\alpha}}\int_{B_1(0)}|v_m(x)|^p\,dx\,dz \geq\frac{\delta}{2}\int_{\mathbb{R}^N}\frac{|v_m(x)|^p}{|x-z|^{N-\alpha}}\,dz,
%\end{eqnarray*}
%combining with
%\begin{eqnarray*}
%\int_{B_1(0)}|v_m(x)|^q\,dx>\frac{\delta}{2},
%\end{eqnarray*}
%we have that
%\begin{eqnarray*}
%\int_{\mathbb{R}^N}\frac{|v_m(x)|^p}{|x-z|^{N-\alpha}}\,dz>0,
%\int_{B_1(0)}\int_{\mathbb{R}^N}\frac{|v_m(x)|^p}{|x-z|^{N-\alpha}}|v_m(x)|^p\,dzdx>0,
%\end{eqnarray*}
%then   $\int_{\mathbb{R}^N}(I_\alpha*|w|^p)(x)|w(x)|^p\,dx\neq 0$.
\medskip
Using  Lemma \ref{lemma 2}, we obtain that
\begin{eqnarray*}
&&\int_{\mathbb{R}^N}(I_\alpha*|w_m|^p)(x)|w_m(x)|^p\,dx
\\&= & \int_{\mathbb{R}^N}(I_\alpha*|w|^p)(x)|w(x)|^p\,dx+\lim_{m\to \infty}\int_{\mathbb{R}^N}(I_\alpha*|w_m-w|^p)(x)|(w_m-w)(x)|^p\,dx
\end{eqnarray*}
and
\begin{eqnarray*}
&& M_p=\lim_{m\to\infty}{\|w_m\|_{H^\frac{1}{2}(\mathbb{R}^N)}^2}={\|w\|_{H^\frac{1}{2}(\mathbb{R}^N)}^2}+\lim_{m\to\infty}{\|w_m-w\|_{H^\frac{1}{2}(\mathbb{R}^N)}^2}
\\&\geq &M_p\left(\lim_{m\to\infty}\int_{\mathbb{R}^N}(I_\alpha*|w_m-w|^p)(x)|(w_m-w)(x)|^p\,dx\right)^\frac{2}{p+1}
\\&& \ \ +M_p\left(\int_{\mathbb{R}^N}(I_\alpha*|w|^p)(x)|w(x)|^p\,dx\right)^\frac{2}{p+1}
\\&=& M_p\left(\int_{\mathbb{R}^N}(I_\alpha*|w|^p)(x)|w(x)|^p\,dx\right)^\frac{2}{p+1}+M_p\left(1-\int_{\mathbb{R}^N}(I_\alpha*|w|^p)(x)|w(x)|^p\,dx\right)^\frac{2}{p+1}
\end{eqnarray*}
Then  $\int_{\mathbb{R}^N}(I_\alpha*|w|^p)(x)|w(x)|^p\,dx=1$. As a  consequent, we get that $M_p=\|w\|_{H^\frac{1}{2}(\mathbb {R}^N)}^2$.\\
The proof is completed. \hfill$\Box$

\setcounter{equation}{0}
\section{Regularity}
To consider the regularity of solutions to problem (\ref{eq 1}),
 we transform (\ref{eq 1}) to the following extension problem
\begin{equation}\label{eq 123}
\left\{ \arraycolsep=1pt
\begin{array}{lll}
  -\Delta w(x,y)+w(x,y)=0,\quad \ (x,y)\in\mathbb{R}_+^{N+1},\\[2mm]
\frac{\partial w}{\partial \nu}(x,0)=(I_\alpha*|w|^p)(x,0)|w(x,0)|^{p-2}w(x,0)\ \  \ \ x\in\mathbb{R}^{N},
\end{array}
\right.
\end{equation}
as in the work of \cite {TWY}.
\begin{proposition}\label{proposition 4.1}
Suppose that $w\in H^1(\mathbb{R}_+^{N+1})$ is a weak solution of  (\ref{eq 123}). Then $w\in L^q_{loc}(\mathbb{R}_+^{N+1})$. Moreover, $w\in C^{2,\alpha}(\mathbb{R}_+^{N+1})$.
\end{proposition}
\noindent {\bf Proof.}
Let $\varphi\in C^\infty(\mathbb{R}^{N+1})$, $0\leq\varphi\leq1$ and for small fixed $R$,
\begin{equation*}
\varphi(x,y) =\left\{ \arraycolsep=1pt
\begin{array}{lll}
 \displaystyle  1  \quad
 \quad  (x,y)\in B_R(0),\\[2mm]
 \phantom{    }
 \displaystyle  0\quad
 \quad  (x,y)\not\in B_{2R}(0),
\end{array}
\right.\qquad
|\nabla \varphi|\leq\frac{C}{R}.
\end{equation*}
Multiplying (\ref{eq 123}) by $\varphi^2|w|^{2\beta_0}w$, where $\beta_0>0$ ,and integrating by parts, we have that
\begin{equation}\label{eq 1234}
\int_{\mathbb{R}_+^{N+1}}(\nabla w\nabla(\varphi^2|w|^{2\beta_0}w)+\varphi^2|w|^{2\beta_0+2})\,dxdy=\int_{\mathbb{R}^N}(I_\alpha*|w|^p)(x,0)\varphi^2(x,0)|w(x,0)|^{2\beta_0+p}\,dx.
\end{equation}
Since
\begin{eqnarray*}
&&\int_{\mathbb{R}_+^{N+1}}\nabla w\nabla(\varphi^2|w|^{2\beta_0}w)\,dxdy
\\&=&\int_{\mathbb{R}_+^{N+1}}\left((2\beta_0+1))\varphi^2|w|^{2\beta_0}|\nabla w|^2+2\varphi|w|^{2\beta_0}w\nabla \varphi\nabla w\right)\,dxdy.
\end{eqnarray*}
So we have that
\begin{eqnarray*}
&&\int_{\mathbb{R}_+^{N+1}}(2\beta_0+1))\varphi^2|w|^{2\beta_0}|\nabla w|^2\,dxdy
\\&&=\int_{\mathbb{R}^N}(I_\alpha*|w|^p)\varphi^2|w|^{2\beta_0+p}\,dx-\int_{\mathbb{R}_+^{N+1}}\left(\varphi^2|w|^{2\beta_0+2}+2\varphi|w|^{2\beta_0}w\nabla \varphi\nabla w\right)\,dxdy
\\&&\leq\int_{\mathbb{R}^N}(I_\alpha*|w|^p)\varphi^2|w|^{2\beta_0+p}\,dx+\int_{\mathbb{R}_+^{N+1}}2\varphi|w|^{2\beta_0+1}w\nabla \varphi\nabla w\,dxdy
\\&&\leq\int_{\mathbb{R}^N}(I_\alpha*|w|^p)\varphi^2|w|^{2\beta_0+p}\,dx+ \left(\int_{\mathbb{R}_+^{N+1}}2\varepsilon\varphi^2|w|^{2\beta_0}w|\nabla \varphi|^2+2C_\varepsilon w^{2(\beta_0+1)}|\nabla \varphi|^2\right)\,dxdy.
\end{eqnarray*}
Note that
\begin{eqnarray}
&&\int_{\mathbb{R}_+^{N+1}}|\nabla(\varphi|w|^{\beta_0}w)|^2\,dxdy\nonumber
\\&\leq&C(\beta_0+1)\int_{\mathbb{R}_+^{N+1}}\left((2\beta_0+1)\varphi^2|w|^{2\beta_0}|\nabla w|^2+2|w|^{2\beta_0+2}|\nabla \varphi|^2\right)\,dxdy,\label{a4}
\end{eqnarray}
so we deduce from (\ref{eq 1234}) and (\ref{a4}) that
\begin{eqnarray*}
&&\int_{\mathbb{R}_+^{N+1}}(|\nabla(\varphi|w|^{\beta_0}w)|^2+\varphi^2|w|^{2\beta_0+2})\,dxdy
\\&&\leq C(\beta_0+1)\int_{\mathbb{R}_+^{N+1}}\left((2\beta_0+1)\varphi^2|w|^{2\beta_0}|\nabla w|^2+2|w|^{2\beta_0+2|\nabla \varphi|^2}+\varphi^2|w|^{2\beta_0+2}\right)\,dxdy
\\&&\leq C(\beta_0+1)\int_{\mathbb{R}^N}(I_\alpha*|w|^p)\varphi^2|w|^{2\beta_0+p}\,dx+ C(\beta_0+1)\int_{\mathbb{R}_+^{N+1}}w^{2(\beta_0+1)}|\nabla \varphi|^2\,dxdy.
\end{eqnarray*}
By the Sobolev theorem, we have that
\begin{eqnarray*}
&&C(R)\left(\int_{\mathbb{R}_+^{N+1}}(|\varphi||w|^{\beta_0+1})^\frac{2N}{N-1}\,dx\right)^\frac{N-1}{N}
\\&&\leq C\left(\int_{\mathbb{R}_+^{N+1}}(|\varphi||w|^{\beta_0+1})^\frac{2N}{N-2}\,dx\right)^\frac{N-2}{N}
\\&&\leq C\int_{\mathbb{R}_+^{N+1}}\left(|\nabla(\varphi|w|^{\beta_0}w)|^2+\varphi^2|w|^{2\beta_0+2})\,dxdy\right)
\\&&\leq C(\beta_0+1)\int_{\mathbb{R}^N}(I_\alpha*|w|^p)\varphi^2|w|^{2\beta_0+p}\,dx+ C(\beta_0+1)\int_{\mathbb{R}_+^{N+1}}w^{2(\beta_0+1)}|\nabla \varphi|^2\,dxdy.
\end{eqnarray*}
By the Hardy-Littlewood-Sobolev inequality, we have that
\begin{eqnarray*}
&&\int_{\mathbb{R}^N}(I_\alpha*|w|^p)\varphi^2|w|^{2\beta_0+p}\,dx
\\&=&\int_{\mathbb{R}^N}\int_{\mathbb{R}^N}\frac{w^p(z)w^{2\beta_0+p}(x)\varphi^2(x)}{|x-y|^{N-\alpha}}\,dzdx\leq\||w|^p\|_{L^{t_1}}\||w|^{2\beta_0+p}\varphi^2\|_{L_{t_2}},
\end{eqnarray*}
where $\frac{1}{t_1}+\frac{1}{t_2}+\frac{N-\alpha}{N}=2$, we let $pt_1=2^\sharp=\frac{2N}{N-1}$ and $\frac{1}{t_2}=1+\frac{\alpha}{N}-\frac{1}{t_1}$.
Note that
$$(2\beta_0+p)t_2\leq \frac{2N}{N-1},$$
we have that
\begin{eqnarray*}
\beta_0\leq\frac{\frac{2N}{N-1}-pt_2}{2t_2}=\frac{N}{N-1}(1+\frac{\alpha}{N}-\frac{p(N-1)}{2N})-\frac{p}{2}=\frac{N+\alpha}{N-1}-p
\end{eqnarray*}
Let $\beta_0=\frac{N+\alpha}{N-1}-p>0$, \ $(2\beta_1+p)t_2=2^\sharp(\beta_0+1)$,\ \ we have that $\beta_1=\frac{2^\sharp(\beta_0+1)-pt_2}{2t_2}>\beta_0\leq\frac{2^\sharp-pt_2}{2t_2}$,
so $\beta_1>\beta_0$.
Let $(2\beta_i+p)t_2=2^\sharp(\beta_{i-1}+1)$.\ Repeating the procedure and using (\ref{eq 1234}), we find $w(x,y)\in L_{loc}^q(\mathbb{R}_+^{N+1})$, where $q\in[2,\infty)$ and
\begin{eqnarray}\label{1234567}
\|w(x,y)\|_{L^q(B_R(0))}\leq C(q,R)\|w(x,y)\|_{L^2(B_{2R}(0))}.
\end{eqnarray}
Now, we use the following auxiliary function
\begin{equation*}
f(x,y)=\int_0^y w(x,t)\,dt\ \ for \ (x,y)\in\mathbb{R}_+^{N+1}.
\end{equation*}
Since $(-\triangle f+f)_y=0$\ \ {\rm in} $\mathbb{R}_+^{N+1}$,\ we have that $-\triangle f+f$ is independent of $y$.\ Hence, if $y=0$,\ $f\equiv 0$, we have
$-\triangle f+f=f_yy=w_y$.\ Thus $w$ is a solution of the Dirichlet problem

\begin{equation}\label{1234567}
\left\{ \arraycolsep=1pt
\begin{array}{lll}
 \displaystyle  -\triangle f(x,y)+f(x,y)=(I_\alpha*|w|^p)(x,0)|w(x,0)|^{p-2}w(x,0)  \quad
 & {\rm for}\quad x\in \mathbb{R}^N,\ y>0\\[2mm]
 \phantom{    }
 \displaystyle  f(x,0)=0\quad
 & {\rm for}\quad  x\in \mathbb{R}^N.
\end{array}
\right.
\end{equation}
Then, by the estimates of Calderon-Zymund, we have $f\in W^{2,p}(\mathbb{R}_+^{N+1})$,\ $w\in W^{1,p}(\mathbb{R}_+^{N+1})$,\ thus,\ the Schauder estimates give
$w\in C^{2,\alpha}(\mathbb{R}_+^{N+1})$.
This proof is complete.\hfill$\Box$

\setcounter{equation}{0}
\section{Nonexistence}
In this section, we prove a Poho\v{z}aev type identity for problem (\ref{eq 123}), which implies the non-existence results in Theorem \ref{th1}.
\begin{lemma}\label{lemma 6}
Let $w\in H^1(\mathbb{R}_+^{N+1})\cap L^\frac{2Np}{N+\alpha}(\mathbb{R}^N)$ be a solution of (\ref{eq 123}) and $\nabla w\in H_{loc}^1(\mathbb {R}^N)$. Then, there holds
\begin{equation}\label{12345678}
\frac{N-1}{2}\int_{\mathbb{R}_+^{N+1}}|\nabla w|^2\,dxdy+\frac{N+1}{2}\int_{\mathbb{R}_+^{N+1}}w^2\,dxdy=\frac{N+\alpha}{2p}\int_{\mathbb{R}^N\times\{0\}}(I_\alpha*|w|^p)|w|^p\,dx.
\end{equation}
\end{lemma}
\noindent {\bf Proof.}
We take $\psi\in C_c^1(\mathbb{R}^N\times[0,\infty))$ such that $\psi=1$ on $B$, where set $B=\{z=(x,y)\in \mathbb{R}^N\times[0,\infty):|z|\leq 1\}$
Let $w\in H^1(\mathbb{R}_+^{N+1})$ be a bounded solutions of  (\ref{eq 123}). By Proposition \ref{proposition 4.1} we know that $w\in C^2(\mathbb{R}_+^{N+1})\cap C^1(\bar{\mathbb{R}_+^{N+1}})$.\ Let $\omega_R=B_R(0)\cap\mathbb{R}_+^{N+1}$, where $B_R(0)\subset\mathbb{R}_+^{N+1}$ is a ball centered at the originn
with the radius $R$.
By (\ref{eq 123}), we have
\begin{equation}\label{12345679}
\int_{\omega_R}(|\nabla w|^2+w^2)\,dxdy =\int_{\partial {w_R}}w\frac{\partial w}{\partial\nu}\,dS.
\end{equation}
Multiplying (\ref{eq 123}) by $(\nabla w, z)$, where $z=(x,y)\in\mathbb{R}_+^{N+1}$, and integrating on $\omega_R$, we deduce that
\begin{eqnarray*}
&&0=\int_{\omega_R}(-\Delta w+w)(\nabla w,z)\,dxdy
\\&=&\int_{\omega_R}\nabla(\frac{1}{2}|\nabla w|)z\,dxdy+\int_{\omega_R}|\nabla w|^2\,dxdy+\frac{1}{2}\int_{\omega_R}(\nabla w^2,z)\,dxdy-\int_{\partial {w_R}}(\nabla w,z)\frac{\partial w}{\partial\nu}\,dS
\\&=&-\frac{N-1}{2}\int_{\omega_R}|\nabla w|^2\,dxdy+\frac{1}{2}\int_{\partial {w_R}}|\nabla w|^2(z,\nu)\,dS-\frac{N+1}{2}\int_{\omega_R}w^2\,dxdy
\\&&+\frac{1}{2}\int_{\partial {w_R}}w^2(z,\nu)\,dS-\int_{\partial {w_R}}(\nabla w,z)\frac{\partial w}{\partial \nu}\,dS.
\end{eqnarray*}
So we have that
\begin{eqnarray}
&&\int_{\omega_R}\frac{N-1}{2}\int_{\omega_R}|\nabla w|^2\,dxdy+\frac{N+1}{2}\int_{\omega_R}w^2\,dxdy\nonumber
\\&=&\frac{1}{2}\int_{\partial {w_R}\cap\{y>0\}}|\nabla w|^2(z,\nu)\,dS+\frac{1}{2}\int_{\partial {w_R}\cap\{y>0\}}w^2(z,\nu)\,dS\nonumber
\\&-&\int_{\partial {w_R}\cap\{y>0\}}(\nabla w,z)\frac{\partial w}{\partial \nu}\,dS-\int_{B_R(0)\times\{0\}}(\nabla_xw,x)\frac{\partial w}{\partial \nu}\,dx\label{a1}
\end{eqnarray}

By the same arguments in the proof of Proposition 3.1 in \cite {VJ}, we obtain that
$$\int_{B_R(0)\times\{0\}}(I_\alpha*|w|^p)|w|^{p-2}w(\nabla_xw,x)\,dx\rightarrow-\frac{N+\alpha}{2p}\int_{B_R(0)\times\{0\}}(I_\alpha*|w|^p)|w|^p\,dx$$
as $R\to\infty$.
Since
$$\left|\int_{\partial {w_R}\cap\{y>0\}}|\nabla w|^2(z,\nu)\,dS\right|\leq CR\int_{\partial {w_R}\cap\{y>0\}}|\nabla w|^2\,dS,$$
$$\left|\int_{\partial {w_R}\cap\{y>0\}}w^2(z,\nu)\,dS\right|\leq CR\int_{\partial {w_R}\cap\{y>0\}}w^2\,dS,$$
we next claim that
$$R\int_{\partial {w_R}\cap\{y>0\}}|\nabla w|^2\,dS\to 0,\ \ \ R\int_{\partial {w_R}\cap\{y>0\}}w^2\,dS\to 0 \quad {\rm as}\quad R\to+\infty.$$
Here we only show that
$$R\int_{\partial {w_R}\cap\{y>0\}}w^2\,dS\to 0\quad {\rm as}\quad R\to+\infty,$$
 the other can be treated in the same way. To this end, by contradiction, we assume that
$$\lim_{R\rightarrow\infty}\inf\int_{\partial {w_R}\cap\{y>0\}}w^2\,dS=c_0>0,$$
then there exists $R_0$ such that, for all $R_1\geq R_0$,
\begin{eqnarray*}
C\geq\int_{\mathbb{R}_+^{N+1}} w^2\,dx\geq
\int_{R_0}^{R_1}\int_{\partial {w_R}\cap\{y>0\}}w^2\,dxdR\geq\frac{c_0}{2}\int_{R_0}^{R_1}\frac{1}{R}\,dR=\frac{c_0}{2}\log\frac{R_1}{R_0}.
\end{eqnarray*}
It yields a contradiction when $R_1>0$ large. So the claim holds.\hfill$\Box$

\medskip

We now complete the proof of  Theorem \ref{th1}$(ii)$.

\medskip

\noindent{\bf Proof of Theorem \ref{th1}$(ii)$.}
Let $w\in H^1(\mathbb{R}_+^{N+1})$ be a solution of  (\ref{eq 123}), we obtain the identity
$$\int_{\mathbb{R}^N}|\nabla w|^2+\int_{\mathbb{R}^N}|w|^2 =\int_{\mathbb{R}^N}(I_\alpha \ast |w|^p) |w|^p .$$
Hence, combine with equation (\ref{eq 123}), we have
$$(\frac{N-1}{2}-\frac{N+\alpha}{2p})\int_{\mathbb{R}^N}|\nabla w|^2\,dx+(\frac{N+1}{2}-\frac{N+\alpha}{2p})\int_{\mathbb{R}^N}w^2\,dx=0.$$
If $0<p\leq\frac{N+\alpha}{N+1}$ or $p\geq\frac{N+\alpha}{N-1}$, then $u = 0$. \hfill$\Box$

\setcounter{equation}{0}
\section{Berestycki-Lions type solutions}

The aim of this section is to establish the infinitely many bounded solutions in Theorem \ref{th2}. We will apply the genus theory to an even functional,
 which is constrained in a manifold and obtain the infinitely many critical points of the functional. To this end,
 let us recall the  following  critical point theorems in  \cite {HP}.

Let $H$ be a real Hilbert space whose norm and inner product will be denoted respectively by $\|\cdot\|$ and $(\cdot,\cdot)$. Consider the manifold
$$\mathcal{M}:=\{w\in H:\|w\|_H=1\},$$
the tangent space of $\mathcal{M}$ at a given point $g\in \mathcal{M}$ is given by
$$T_w\mathcal{M}=\{w\in H:(g,w)=0\}.$$
Let $J$ be a $C^1$ functional defined on $H$. Then the trace $J|_\mathcal{M}$ of $J$ on $\mathcal{M}$ is of class $C^1$ and for any $w\in \mathcal{M}$,
$$\langle J|_\mathcal{M}'(w),g\rangle=\langle J'(w),g\rangle \ \ \  \ g\in T|_w\mathcal{M}.$$
Let $\sum(\mathcal{M})$ denote the set of compact and symmetric subsets of $\mathcal{M}$ . The genus $\gamma(A)$ of a set $A\in \sum(\mathcal{M})$ is defined as the least integer $n\geq1$ such that there exists an odd continuous mapping $\varphi:A\rightarrow S^{n-1}$. For $k\geq1$, we denote $\tau_k=\{A\in \sum(\mathcal{M}):\gamma(A)\geq k\}$.

We say that a functional $J$ defined on a manifold $\mathcal{M}$ satisfies the positive Palais-Smale condition (in short,$(PS)^+$) if for $0<c_1<c_2$, for every sequence $\{w_m\}\subset\mathcal{M}$ such that $\|J|_\mathcal{M}'(w_m)\|\rightarrow 0$ and $c_1\leq J(w_m)\leq c_2$, there exists a convergent subsequence $\{w_{m_i}\}$ of $\{w_m\}$.
\begin{proposition}\label{proposition 6.1}
Let $J:H\rightarrow R$ be an even functional of class $C^1$. Suppose that J is bounded from above on $\mathcal{M}$  and $J|_\mathcal{M}$
satisfies the (PS) condition. Let
$$b_k=\sup_{A\in\tau_k}\inf_{w\in A}J(w).$$
Then for any $k\geq1$, $b_k$ is a critical value of $J|_\mathcal{M}$ and $b_1\geq b_2\geq\cdots \geq b_k\geq\cdots$. If J only satisfies the $(PS)^+$ condition, then $b_k$
is a critical value of provided $b_k>0$.
\end{proposition}
To be convinent for the analysis, we denote
$$ K_p=\left\{x\in \mathcal{M}: J(w)=b,J|_\mathcal{M}'(w)=0\right\}.$$
\begin{proposition}\label{proposition 6.2}
 Under the hypotheses of Proposition \ref{proposition 6.1}, suppose that $ b_k=b_{k+1}=\cdots=b_{k+r-1}\equiv b$. Then $\gamma(K_b)\geq r$. In particular, if $r\geq2$,
 there exist infinitely many distinct critical points of $J|_\mathcal{M}$ corresponding to the critical value $b$.
\end{proposition}
By Proposition \ref{proposition 6.1} and Proposition \ref{proposition 6.2}, under the conditions of Proposition \ref{proposition 4.1},
there always exist infinitely many distinct critical points of $J $  on $\mathcal{M}$.
Let $E=H^\frac{1}{2}(\mathbb{R}^N)$ and $E_r=\left\{v\in H^\frac{1}{2}(\mathbb{R}^N): \ v(x)=v(|x|)\right\}.$
 Denote by $\mathcal{M}=\left\{v\in E_r: \ \|v\|_{E_r}=1\right\}$ the unit ball in $E_r$. Define the functional
 $$J(v)=\int_{\mathbb{R}^N}(I_\alpha*|v|^p)(x)|v(x)|^p\,dx\ \ \  for \ v\in E_r.$$
 We next verify that $J$  satisfies the conditions of Proposition \ref{proposition 6.1}. In fact, by Hardy-Littlewood-Sobolev inequality, we have that
 \begin{eqnarray*}
&&\int_{\mathbb{R}^N}(I_\alpha*|v|^p)(x)|v(x)|^p\,dx=\int_{\mathbb{R}^N}\int_{\mathbb{R}^N}\frac{|v(z)|^p|v(x)|^p}{|x-z|^{N-\alpha}}\,dxdz
\\&& \leq c\||v|^p\|_l\|v|^p\|_t=c\left(\int_{\mathbb{R}^N}|v(x)|^\frac{2Np}{N+\alpha}\,dx\right)^\frac{N+\alpha}{N},
\end{eqnarray*}
where $l=t=\frac{2N}{N+\alpha}$.
 Then, by Min-Max method argument, we have that
 \begin{eqnarray*}
&&\int_{\mathbb{R}^N}(I_\alpha*|v|^p)(x)|v(x)|^p\,dx
\\&\leq &c_6\left(\sup_{z\in{\mathbb{R}^N}}(\int_{B_r(z)}|v(x)|^q\,dx)^{\frac{1}{q}(\frac{2Np}{N+\alpha}-2)}\int_{\mathbb{R}^N}
|(-\Delta+id)^\frac{1}{4}u|^2\,dx\right)^\frac{N+\alpha}{N},
\end{eqnarray*}
combining with the fact that $\|v\|_{E_r}=1$ and $H^\frac{1}{2}(\mathbb{R}^N) \hookrightarrow L^p(\mathbb{R}^N)$ for $2\leq p\leq \frac{2N}{N-1}$,
we have that $\int_{B_r(z)}|v(x)|^q\,dx<+\infty$ and
\begin{eqnarray*}
\sup_{z\in{\mathbb{R}^N}}\int_{B_r(z)}|v(x)|^q\,dx<+\infty .
\end{eqnarray*}
Thus, $J$ is bounded from above on $\mathcal{M}$.
\begin{lemma}\label{lemma 4.1}
Let $2<q<2^\sharp:=\frac{2N}{N-1}$ for $N\geq 2$, $\{\tilde{v}_m\}_m$ be a bounded sequence in $E_r$,
suppose that $\tilde{v}_m\rightharpoonup 0$ in $H^\frac{1}{2}(\mathbb{R}^N)$ as $m\to\infty$. Then
$$\int_{\mathbb{R}^N}(I_\alpha*|\tilde{v}_m|^p)(x)|\tilde{v}_m(x)|^p\,dx\to 0$$
as $m\to\infty$.
\end{lemma}
\noindent{\bf Proof}.
We know that
$\int_{\mathbb{R}^N}(I_\alpha*|\tilde{v}_m|^p)(x)|\tilde{v}_m|^p\,dx$ strongly converges, up to a subsequence.
Denote by $\{v_m\}$ the restriction of
functions  $\{\tilde{v}_m\}$ to $\mathbb{R}^N$. Then $\{v_m\}\subset H^\frac{1}{2}(\mathbb{R}^N)$ is uniformly bounded  and $v_m$ is radially
symmetric in x. Assuming $v_m\rightharpoonup 0$ in $ H^\frac{1}{2}(\mathbb{R}^N)$ as $m\to \infty$.
By the results in \cite {TWY}, we have that $$\sup_{z\in{\mathbb{R}^N}}\int _{B_r(z)}|v_m(x)|^2\,dx\to 0$$ as $m\to\infty$. Combining with Lemma \ref{lemma 2.1}, it implies that
$$\int_{\mathbb{R}^N}(I_\alpha*|v_m|^p)(x)|v_m(x)|^p\,dx \to 0$$ as $m\to\infty$, \ and then $$\int_{\mathbb{R}^N}(I_\alpha*|\tilde{v}_m|^p)(x)|\tilde{v}_m(x)|^p\,dx \to 0$$
as $m\to \infty$ by the fact that $\tilde{v}_m(x)=v_m(x)$.
\begin{lemma}\label{lemma 4.2}
$J|_\mathcal{M}$ satisfies the $(PS)^+$ condition.
\end{lemma}
\noindent{\bf Proof.}
Let $\{v_m\}\subset \mathcal{M}$ be a $(PS)^+$ sequence for $J$, that is,
$$0<\alpha \leq J(v_m)\leq C,\quad J|_\mathcal{M}'(v_m)\to 0$$
 as $m\to \infty$.
Since $v_m$ is bounded in $E_r$, we may assume $v_m\rightharpoonup v$ in $E_r$. Then
$$<J'(v_m)-<J'(v_m),v_m>v_m,\varphi>\to 0$$
 as $m\to\infty$ and for any $\varphi\in E_r$,
\begin{eqnarray*}
\int_{\mathbb{R}^N}(I_\alpha*|v_m|^p)(x)|v_m(x)|^{p-2}v_m(x)\varphi(x)\,dx-(\int_{\mathbb{R}^N}(I_\alpha*|v_m|^p)(x)|v_m(x)|^p\,dx)\cdot\\
\int_{\mathbb{R}^N}(\nabla v_m(x)\nabla\varphi(x)+v_m(x)\varphi(x))\,dx=o(1).
\end{eqnarray*}

By lemma \ref{lemma 4.1}, we have that
$$\int_{\mathbb{R}^N}(I_\alpha*|v_m|^p)(x)|v_m(x)|^p\,dx\to \int_{\mathbb{R}^N}(I_\alpha*|v|^p)(x)|v(x)|^p\,dx \quad {\rm as}\quad m\to\infty,$$
together with the weak convergence of $\{v_m\}$,  yield
$$\int_{\mathbb{R}^N}(\nabla v_m(x)\nabla\varphi(x)+v_m(x)\varphi(x))\,dx\to\int_{\mathbb{R}^N}(\nabla v(x)\nabla\varphi(x)+v(x)\varphi(x))\,dx \quad {\rm as}\quad m\to\infty$$
and
$$\int_{\mathbb{R}^N}(I_\alpha*|v_m|^p)(x)|v_m(x)|^{p-2}v_m(x)\varphi(x)\,dx\to\int_{\mathbb{R}^N}(I_\alpha*|v|^p)(x)|v(x)|^{p-2}v(x,0)\varphi(x)\,dx.$$
Denote that $\xi^m=(I^\alpha*|v_m|^p)(x)|v_m(x)|^{p-1}$, $\xi=(I^\alpha*|v|^p)(x)|v(x)|^{p-1}$,
we have that $\xi^m\rightharpoonup \xi$ in $E_r$ and
\begin{eqnarray*}
&&\int_{\mathbb{R}^N}(I_\alpha*|v_m|^p)(x)|v_m(x)|^{p-1}\varphi(x)\,dx-\int_{\mathbb{R}^N}(I_\alpha*|v|^p)(x)|v(x)|^{p-1}\varphi(x)\,dx
\\&=&\int_{\mathbb{R}^N}(I_\alpha*(|v_m|^p-|v|^p))(x)|v_m(x)|^{p-1}\varphi(x)\,dx+\int_{\mathbb{R}^N}(I_\alpha*|v|^p)(x)(|v_m(x)|^{p-1}-|v(x)|^{p-1})\varphi(x)\,dx
\\&\leq &C_1\int_{\mathbb{R}^N}(I_\alpha*(|v_m|^{p-1}+|v|^{p-1})|v_m-v|)(x)|v_m(x)|^{p-1}\varphi(x)\,dx
\\&&\ \ +C_2\int_{\mathbb{R}^N}(I_\alpha*|v^p|)(x)((|v_m|^{p-2}+|v|^{p-2})|(v_m-v)|)(x)|\varphi(x)\,dx
\\&\leq &C_1\left(\int_{\mathbb{R}^N}(|v_m|^{p-1}+|v|^{p-1})^{\frac{1}{1+\frac{\alpha}{N}-\frac{p}{2^\sharp}}}|(v_m-v)(x)|^
{\frac{1}{1+\frac{\alpha}{N}-\frac{p}{2^\sharp}}}\,dx\right)^{1+\frac{\alpha}{N}-\frac{p}{2^\sharp}}
\left(\int_{\mathbb{R}^N}(|v_m(z)|^{p-1}\varphi(z))^{\frac{2^\sharp}{p}}\,dz\right)^\frac{p}{2^\sharp}
\\&&\ \ +C_2\left(\int_{\mathbb{R}^N}|v^p(z)|^{\frac{2^\sharp}{p}}\,dz\right)^\frac{p}{2^\sharp}\left(\int_{\mathbb{R}^N}
[((|v_m|^{p-2}+|v|^{p-2})|(v_m-v)|)(x)|\varphi(x)]^{\frac{1}{1+\frac{\alpha}{N}-\frac{p}{2^\sharp}}}\,dx\right)^{1+\frac{\alpha}{N}-\frac{p}{2^\sharp}},
\end{eqnarray*}
be the fact that
\begin{eqnarray*}
&&\left(\int_{\mathbb{R}^N}\left((|v_m|^{p-1}+|v|^{p-1})|(v_m-v)(x)|\right)^{\frac{1}{1+\frac{\alpha}{N}-\frac{p}{2^\sharp}}}\,dx\right)^{1+\frac{\alpha}{N}-\frac{p}{2^\sharp}}
\\&\leq &\left(\int_{\mathbb{R}^N}|(v_m-v)(x)|^{\frac{2^\sharp}{(1+\frac{\alpha}{N}-\frac{p}{2^\sharp})2^\sharp-p+1}}\,dx\right)^\frac{(1+\frac{\alpha}{N}-\frac{p}{2^\sharp})2^\sharp-p+1}{2^\sharp}\left(\int_{\mathbb{R}^N}(|v_m|^{p-1}+|v|^{p-1})^\frac{2^\sharp}{p-1}\,dx\right)^\frac{p-1}{2^\sharp}.
\end{eqnarray*}
Since $v_m\rightharpoonup v$ in $E_r$,  we have that $E_r$ is compactly embedded in $L_q(\mathbb{R}^N)$, where $2\leq q< 2^\sharp$ and ${\frac{2^\sharp}{(1+\frac{\alpha}{N}-\frac{p}{2^\sharp})2^\sharp-p+1}}<2^\sharp$ for $p<\frac{N+\alpha}{N-1}$ and $(p-1)\frac{2^\sharp}{p-1}=2^\sharp$,
then
$$\left(\int_{\mathbb{R}^N}(|v_m|^{p-1}+|v|^{p-1})^{\frac{2^\sharp}{p-1}}\,dx\right)^\frac{p-1}{2^\sharp}<\infty,$$
$$\left(\int_{\mathbb{R}^N}|(v_m-v)(x)|^{\frac{2^\sharp}{(1+\frac{\alpha}{N}-\frac{p}{2^\sharp})2^\sharp-p+1}}\,dx\right)^
\frac{(1+\frac{\alpha}{N}-\frac{p}{2^\sharp})2^\sharp-p+1}{2^\sharp}\rightarrow 0$$
and
\begin{eqnarray*}
\left(\int_{\mathbb{R}^N}(|v_m(z)|^{p-1}\varphi(z))^\frac{2^\sharp}{p}\,dz\right)^\frac{p}{2^\sharp}\leq \left(\int_{\mathbb{R}^N}|v_m(z)|^{2^\sharp}\,dz\right)^\frac{p-1}{2^\sharp}\left(\int_{\mathbb{R}^N}\varphi(z)^{2^\sharp}\,dz\right)^\frac{1}{2^\sharp}
\end{eqnarray*}
we  know that $E_r$ is embedded in $L_q(\mathbb{R}^N)$, where $2\leq q\leq2^\sharp$, so $$\left(\int_{\mathbb{R}^N}(|v_m(z)|^{p-1}\varphi(z))^\frac{2^\sharp}{p}\,dz\right)^\frac{p}{2^\sharp}<+\infty,$$
thus, we have that
\begin{eqnarray*}
&&\left(\int_{\mathbb{R}^N}|v^p(z)|^\frac{2^\sharp}{p}\,dz\right)^\frac{p}{2^\sharp}\left(\int_{\mathbb{R}^N}[((|v_m|^{p-2}+|v|^{p-2})|(v_m-v)|)(x)|\varphi(x)]^\frac{1}{1+\frac{\alpha}{N}-\frac{p}{2^\sharp}}\,dx\right)^{1+\frac{\alpha}{N}-\frac{p}{2^\sharp}}
\\&\leq &\left(\int_{\mathbb{R}^N}((|v_m|^{p-2}+|v|^{p-2})(x)\varphi(x))^\frac{2^\sharp}{p-1}\,dx\right)^\frac{p-1}{2^\sharp}\left(\int_{\mathbb{R}^N}|(v_m-v)|)(x)|^{\frac{2^\sharp}{(1+\frac{\alpha}{N}-\frac{p}{2^\sharp})2^\sharp-p+1}}\,dx\right)^\frac{(1+\frac{\alpha}{N}-\frac{p}{2^\sharp})2^\sharp-p+1}{2^\sharp}
\\&\leq &\left(\int_{\mathbb{R}^N}((v_m|^{p-2}+|v|^{p-2})(x))^\frac{2^\sharp}{p-2}\,dx\right)^\frac{p-2}{2^\sharp}
\left(\int_{\mathbb{R}^N}\varphi(x)^{2^\sharp}\,dx\right)^\frac{1}{2^\sharp}\cdot
\\&&\ \ \left(\int_{\mathbb{R}^N}|(v_m-v)|)(x)|^{\frac{2^\sharp}{(1+\frac{\alpha}{N}-\frac{p}{2^\sharp})2^\sharp-p+1}}\,dx\right)^\frac{(1+\frac{\alpha}{N}-\frac{p}{2^\sharp})2^\sharp-p+1}{2^\sharp}
\end{eqnarray*}
where ${\frac{2^\sharp}{(1+\frac{\alpha}{N}-\frac{p}{2^\sharp})2^\sharp-p+1}}<2^\sharp$, \ as $p<\frac{N+\alpha}{N-1}$, \ $(p-2)\frac{2^\sharp}{p-2}=2^\sharp$,\ $p\frac{2^\sharp}{p}=2^\sharp$.
So we have that
$$\left(\int_{\mathbb{R}^N}|(v_m-v)|)(x)|^{\frac{2^\sharp}{(1+\frac{\alpha}{N}-\frac{p}{2^\sharp})2^\sharp-p+1}}\,dx\right)^\frac{(1+\frac{\alpha}{N}-\frac{p}{2^\sharp})
2^\sharp-p+1}{2^\sharp}\to 0,$$ $$\left(\int_{\mathbb{R}^N}((v_m|^{p-2}+|v|^{p-2})(x))^\frac{2^\sharp}{p-2}\,dx\right)^\frac{p-2}{2^\sharp}
\left(\int_{\mathbb{R}^N}\varphi(x)^{2^\sharp}\,dx\right)^\frac{1}{2^\sharp}<+\infty.$$
So we have that
$$\int_{\mathbb{R}^N}(I_\alpha*|v_m|^p)(x)|v_m(x)|^{p-2}v_m(x)\varphi(x)\,dx\to \int_{\mathbb{R}^N}(I_\alpha*|v|^p(x))|v(x)|^{p-2}v(x)\varphi(x)\,dx$$
and
\begin{equation*}
\int_{\mathbb{R}^N}(I_\alpha*|v|^p)(x)|v(x)|^{p-2}v(x)\varphi(x)\,dx-\int_{\mathbb{R}^N}(I_\alpha*|v|^p)(x)|v(x)|^p\,dx
\int_{\mathbb{R}^N}\left(\nabla v\nabla\varphi+v\varphi\right)\,dx=0.
\end{equation*}
Notice that $J(v_m)\geq \alpha>0$, then we have that $J(v_m)\geq\alpha$. So $v\neq0$ and $\|v\|_{E_r}=1$. Combining  $v_m\rightharpoonup v$ with $\|v_m\|_{E_r}\rightarrow\|v\|_{E_r}$,
we obtain that  $v_m$ strongly converges to $v$ in $E_r$.

Next we show that $ b_k>0 $ for each $k\geq1$. Indeed, for $k\geq1$,
we denote
$$\pi_{k-1}=\{l=(l_1,l_2,...l_k)\in \mathbb{R}^k\ \ | \sum_{i=1}^k|l_i|=1\}.$$
Since $\pi_{k-1}$ is homeomorphic to $S^{k-1}$ by an odd homeomorphism, it follows that $\gamma(\pi_{k-1})=k$. The following result is due to Berestycki and Lions.

\begin{lemma}\label{lemma 4.3}
For all $k\geq1$, there exists a constant $R=R(k)>0$ and an odd continuous mapping $\tau: \pi_{k-1}\to H_0^1(B_R(0))$
such that

(i) $\pi_(l)$ is a radial function for all $l\in \pi_{k-1}$ and $0\nsubseteq \tau(\pi_{k-1})$;

(ii) there exist $p,C>0$ such that $\rho\leq\|\nabla u\|_{L^2(B_{R(0)})}\leq C$ for $u\in\tau(\pi_{k-1})$;

(iii) for $u\in\tau(\pi_{k-1}) $, $\int_{\mathbb{R}^N}|u(x)|^q\,dx\geq1$.

\end{lemma}

\begin{lemma}\label{lemma 4.4}
When $N\geq3$, there holds $b_k>0$ for each $k\geq1$.
\end{lemma}
\noindent{\bf Proof.}
Let $\hat{\pi}_{k-1}=\tau(\pi_{k-1})$ and define for each $u\in H_0^1(B_R(0))$ an extension $\tilde{u}\in H^1(\mathbb{R}^N)$ of u such that $\tilde{u}=u$
on $B_R(0)$ and $\tilde{u}=0$ on $\mathbb{R}^N\backslash B_R(0)$.
Denote $\tilde{u}_\sigma(x)=\tilde{u}(\frac{x}{\sigma})$. For $u\in \hat{\pi}_{k-1}$, we define $\chi(u)=(1-y)(I_\alpha*|\tilde{u}_\sigma|^p(x))|\tilde{u}_\sigma(x)|^p$ and $\chi(u)=0$
if $y\not\in R\backslash(0,1)$,
where $\sigma=\sigma(u)>0$ is determined by requiring $\chi(u)\in M$, that is,
$$\int_{\mathbb{R}^N}(|\nabla\chi(u))|^2+|\chi|^2)\,dx=1.$$
We deduce that
\begin{eqnarray*}
&&1=\int_{\mathbb{R}_+^{N+1}}(|\nabla\chi(u))|^2+|\chi|^2)\,dxdy
\\&&=\int_0^1\,dy\int_{\mathbb{R}^N}|(\nabla_x[(I_\alpha*|\tilde{u}_\sigma|^p)(x)|\tilde{u}_\sigma(x)|^p],-(I_\alpha*|\tilde{u}_\sigma|^p)(x)|\tilde{u}_\sigma(x)|^p)|^2\,dx\\
\\&&+\int_0^1\,dy\int_{\mathbb{R}^N}(1-y)^2|(I_\alpha*|\tilde{u}_\sigma|^p)(x)|\tilde{u}_\sigma(x)|^p|^2\,dx
\\&&=\int_0^1\,dy\int_{B_{\sigma R}(0)}|(\nabla_x[(I_\alpha*|\tilde{u}_\sigma|^p)(x)|\tilde{u}_\sigma(x)|^p],
-(I_\alpha*|\tilde{u}_\sigma(x)|^p)|\tilde{u}_\sigma(x)|^p)|^2\,dx\\
\\&&+\int_0^1\,dy\int_{B_{\sigma R}(0)}(1-y)^2|(I_\alpha*|\tilde{u}_\sigma|^p)(x)|\tilde{u}_\sigma(x)|^p|^2\,dx
\end{eqnarray*}
and
\begin{eqnarray*}
\nabla_x[(I_\alpha*|\tilde{u}_\sigma|^p)(x)|\tilde{u}_\sigma(x)|^p]
&=&\nabla_x[\int_{\mathbb{R}^N}\frac{|u(z)|^p}{|\frac{x}{\sigma}-z|^{N-\alpha}}\,dz|u(\frac{x}{\sigma})|^p]
\\&=&\int_{\mathbb{R}^N}|u(z)|^p(\alpha-N)|\frac{x}{\sigma}-z|^{\alpha-N-2}(\frac{x}{\sigma}-z)\frac{1}{\sigma}\,dz|u(\frac{x}{\sigma})|^p
\\&&+\int_{\mathbb{R}^N}\frac{|u(z)|^p}{|\frac{x}{\sigma}-z|^{N-\alpha}}\,dz p|u(\frac{x}{\sigma})|^{p-2}u(\frac{x}{\sigma})\frac{1}{\sigma}
\end{eqnarray*}
and
\begin{eqnarray*}
&&\int_{B_{\sigma R}(0)}|\int_{\mathbb{R}^N}|u(z)|^p(\alpha-N)|\frac{x}{\sigma}-z|^{\alpha-N-2}(\frac{x}{\sigma}-z)\frac{1}{\sigma}\,dz|u(\frac{x}{\sigma})|^p
\\&&+\int_{\mathbb{R}^N}\frac{|u(z,0)|^p}{|\frac{x}{\sigma}-z|^{N-\alpha}}\,dz p|u(\frac{x}{\sigma})|^{p-2}u(\frac{x}{\sigma})\frac{1}{\sigma}|^2\,dx
\\&=&\int_{B_R(0)}|\int_{\mathbb{R}^N}|u(z)|^p(\alpha-N)|y-z|^{\alpha-N-2}(y-z)\frac{1}{\sigma}\,dz|u(y)|^p
\\&&+\int_{\mathbb{R}^N}\frac{|u(z,0)|^p}{|y-z|^{N-\alpha}}\,dz p|u(y)|^{p-2}u(y)\frac{1}{\sigma}|^2\sigma^N\,dy
\\&=&\frac{1}{\sigma^2}\sigma^N\int_{B_R(0)}|\nabla_x[(I_\alpha*|u|^p(x))|u(x)|^p]|^2\,dx.
\end{eqnarray*}
Then
\begin{eqnarray*}
&&\int_{\mathbb{R}_+^{N+1}}(|\nabla\chi(u))|^2+|\chi|^2)\,dxdy=\sigma^{N-2}\int_{B_R(0)}|\nabla_x[(I_\alpha*|u|^p(x,0))|u(x,0)|^p]|^2\,dx
\\&&+\frac{4}{3}\sigma^N\int_{B_R(0)}|(I_\alpha*|u|^p(x,0))|u(x,0)|^p|^2\,dx.
\end{eqnarray*}
Since the right hand side is increasing in $\sigma>0$, we find a unique $\sigma>0$ so that $\chi\in M$. By (ii) of lemma \ref{lemma 4.3}, we obtain
$$1\geq\sigma^{N-2}\int_{B_R(0)}|\nabla_x u|^2\,dx\geq\sigma^{N-2}\rho^2.$$
Hence, there exists $\bar{\sigma}> 0$ independent of $u\in\hat{\pi}_{k-1}$, such that $\sigma(u)\leq\bar{\sigma}$.\ By Poincare's inequality, there holds
$$1\leq C(\sigma^{N-2}+\sigma^N),$$\\
which implies that $\sigma(u)$ has a lower bound $\bar{\sigma}>0$ independent of u.

 Now we prove $b_k>0$ for each $k\geq1$. In fact, we observe that
\begin{eqnarray*}
\int_{\mathbb{R}^N\times\{0\}}|\chi(u)|^p\,dx
&=&\int_{\mathbb{R}^N\times\{0\}}|(1-y)(I_\alpha*|\tilde{u}_\sigma|^p)(x,0)|\tilde{u}_\sigma(x,0)|^p|^p\,dx
\\&=&\int_{\mathbb{R}^N\times\{0\}}|(1-y)|^p|(I_\alpha*|\tilde{u}_\sigma|^p)(x,0)|\tilde{u}_\sigma(x,0)|^p|^p\,dx
\\&=&\int_{\mathbb{R}^N}|(I_\alpha*|\tilde{u}_\sigma|^p)(x,0)|\tilde{u}_\sigma(x,0)|^p|^p\,dx
\\&=&\int_{\mathbb{B}_{\sigma R}}|(I_\alpha*|u(\frac{x}{\sigma},0)|^p)(x,0)|u(\frac{x}{\sigma},0)|^p|^p\,dx
\\&=&\int_{\mathbb{B}_R}|(I_\alpha*|u(y,0)|^p)(x,0)|u(y,0)|^p|^p\,dy
\\&\geq&\sigma^N\cdot1
\ge\underline{\sigma}^N.
\end{eqnarray*}
Let $A_k=(1-y)\chi(\hat{\pi}_{k-1})$ for $0\leq y\leq1$. Since $\chi$ is an odd continuous mapping, $A_k\in \Sigma(M)$.
 Furthermore, $(1-y)\chi\cdot\tau:\pi_{k-1}\rightarrow A_k$
is odd and continuous, so we have that $\gamma(A_k)\geq\gamma(\pi_{k-1})=k$. Hence, $A_k\in\Gamma_k$. Set $\varepsilon_k=\underline{\sigma}^N$, we have that
$$\int_{\mathbb{R}^N\times\{0\}}(I_\alpha*|u|^p)|u|^{p-2}u\,dx\geq\varepsilon_k,\ \ \ \ u\in A_k.$$
Therefore, for each $k\geq1$,
$$b_k=\sup_{A\in \Gamma_k}\inf_{u\in A_k}J(u)\geq\inf_{u\in A_k}J(u)\geq\varepsilon_k>0.$$

\noindent{\bf Proof of Theorem \ref{th2}.}
By lemma \ref{lemma 4.2},
$J|_\mathcal{M}$ satisfies the $(PS)^+$ condition, and Lemma \ref{lemma 4.4} yields $b_k>0$ for
each $k>1$. Then the conclusion  follows  Proposition  \ref{proposition 6.1} and  Proposition \ref{proposition 6.2}. \hfill$\Box$

\bigskip

\noindent{\bf Acknowledgements:} The author would like to express the  warmest gratitude to Prof. Jianfu Yang,
for proposing the problem and for its active participation.
This work is supported by the Jiangxi Provincial Natural Science Foundation (20161ACB20007).

\end{document}